\documentclass{article}
\usepackage{hyperref}
\usepackage{geometry}
\usepackage{amsfonts}
\usepackage{array}
\usepackage{amsmath,mathtools}
\usepackage{amsthm}
\usepackage[utf8]{inputenc}
\usepackage[english]{babel}
\usepackage{amssymb}
\usepackage{systeme}
\usepackage{tikz}
\usepackage{fancyhdr}

\usepackage{lineno}

\theoremstyle{definition}
\newtheorem{defn}{Definition}
\newtheorem{prop}{Proposition}
\newtheorem{thm}{Theorem}

\newtheorem{prob}{Problem}

\DeclareMathOperator{\cl}{cl}

\usepackage[font=small,labelfont=bf]{caption}
\usepackage[font=small,labelfont=normalfont,labelformat=simple]{subcaption}

\usepackage{xcolor}

\begin{document}
	
	\title{Independent Hyperplanes in Oriented Paving Matroids}
	
	\author{Lamar Chidiac \and Winfried Hochstättler}
	
	\date{}

	\maketitle

	\begin{center}
		{\footnotesize
			Fakult\"at f\"ur Mathematik und Informatik, \\
			FernUniversit\"at in Hagen, Germany,\\
			\text{\{lamar.chidiac,winfried.hochstaettler\}@fernuni-hagen.de}
			\\\ \\
		}
	\end{center}

	\begin{abstract}
		In 1993, Csima and Sawyer \cite{CsimaSawyer1993} proved that in a non-pencil arrangement of n pseudolines, there are at least $\frac{6}{13}n$ simple points of intersection. Since pseudoline arrangements are the topological representations of reorientation classes of oriented matroids of rank $3$, in this paper, we will use this result to prove by induction that an oriented paving matroid of rank $r \ge 3$ on $n$ elements, where $n \geq 5+ r$, has at least $\frac{12}{13(r-1)} \binom{n}{r-2}$ independent hyperplanes, yielding a new necessary condition for a paving matroid to be orientable.
	\end{abstract}

	\section{Introduction}
	\label{sec:intro}
	In 1893 Sylvester asked in \cite{sylvester} the following question:
	``Prove that it is not possible to arrange any finite number of real points so that a right line through every two of them shall pass through a third, unless they all lie in the same right line."
	This question remained unsolved for almost 40 years, until it was independently raised by Erd\"os and solved shortly after by Gallai in 1933 \cite{erdos1,erdos2,erdos3} and later on several other proofs have been found (for reference \cite{motzkin}, p.451 and \cite{coxeter}, p. 65). 
	
	In other words, the Sylvester-Gallai theorem states that given a set of non-collinear points in the Euclidean plane, we can always find at least one line that has exactly two of the given points, we call it an {\em ordinary line}. A generalization of this theorem to higher dimension is not always true, i.e. given a finite set of points in a d-dimensional Euclidean space which is not contained in a hyperplane, we cannot always find a hyperplane containing exactly d of the given points, which we call an {\em independent hyperplane}. A counterexample would be a set of 6 points, 3 on each of two skew lines as seen in Figure~\ref{fig:hansen} below; there is no hyperplane spanned by these points that has exactly 3 of them. This counterexample is by Hansen \cite{hansen65} and another one is the Desargues configuration in 3-space \cite[p.~452]{motzkin}.

	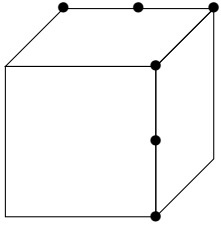
\begin{figure}
		
		\centering
		\begin{tikzpicture}[scale=2.0]
			\pgfmathsetmacro{\cubex}{1}
			\pgfmathsetmacro{\cubey}{1}
			\pgfmathsetmacro{\cubez}{1}
			\draw[black] (0,0,0) -- ++(-\cubex,0,0) -- ++(0,-\cubey,0) -- ++(\cubex,0,0) -- cycle;
			\draw[black] (0,0,0) -- ++(0,0,-\cubez) -- ++(0,-\cubey,0) -- ++(0,0,\cubez) -- cycle;
			\draw[black] (0,0,0) -- ++(-\cubex,0,0) -- ++(0,0,-\cubez) -- ++(\cubex,0,0) -- cycle;
			\foreach \Point in {(0,0,0),(0,-1,0),(0,-1/2,0),(0,0,-1),(-1/2,0,-1),(-1,0,-1)}{
				\node at \Point {\textbullet};
			}
		\end{tikzpicture}
		
		\caption{An illustration of Hansen's construction \cite{hansen65} }
		\label{fig:hansen}
	\end{figure}

	Looking into Hansen's counterexample we notice that the main issue lies in the 3-point lines, and we can actually forbid this by considering a more specific type of point configuration. We consider oriented paving matroids which we will define later. In fact, realizable simple oriented paving matroids of rank 4 are exactly the point configurations in three space that have no three point lines. By considering this type of matroids we can go further with the generalization and not only prove the existence of an independent hyperplane but also examine how many do we have. To do so we will have to go back and start with the 2-dimensional space.
	
	Sylvester-Gallai's theorem can surely be dualized to the projective plane, where the theorem will be equivalent to the following: For every finite set of lines, not all going through one point (non-pencil arrangement), then among all the intersection points of these lines, at least one is incident with exactly two of the lines. We call it a {\em simple point}.
	
	After proving the existence of at least one simple point, Dirac turned to the question of how many. He \cite{dirac1951} proved that there are at least 3 simple points in an arrangements of more than 2 lines and conjectured that there are at least n/2 simple points in an arrangement of n lines. This is a long-standing conjecture which has been known as the {\em Dirac-Motzkin conjecture}. Apparently neither author formally conjectures this in print, but Dirac \cite{dirac1951} states twice that its truth is "likely". Motzkin \cite{motzkin} does not seem to mention it at all. 
	This conjecture was and still is a motivation for a lot of mathematicians to work on a lower bound for the number of simple points in line arrangements but very few results are known since then, the latest was Tao and Green's result in 2013 \cite{taogreen} where they proved Dirac's conjecture when n is sufficiently large. The conjecture however remains open in general.
	
	A careful consideration of the proofs of most of the results for line arrangements reveals that the straightening of the lines which form the arrangements plays only a very limited role. This leads naturally to the idea of investigating arrangements of more general types; arrangements of pseudolines.
	
	\begin{defn}
		An {\em arrangement of pseudolines} is a finite collection of simple curves (no self intersection) in the real projective plane satisfying the following two properties: 
		\begin{enumerate}
			\item any two curves intersect in exactly one point, where they cross.
			\item the intersection of all curves is empty.  
		\end{enumerate}
		As in line arrangements, a simple point in a pseudoline arrangement is the point formed by the intersection of exactly 2 pseudolines (curves). In \cite{CsimaSawyer1993}, simple points are referred to as ordinary points.
	\end{defn}
	
	In the following, we present our main result that concerns a lower bound for the number of independent hyperplanes in an oriented paving matroid. In Section~\ref{sec:application}, we show an application to this result and then give its proof in Section~\ref{sec:proof}. Before going any further, let's see the main connection between oriented matroids and pseudoline arrangements and the main motivation behind paving matroids. 
	
	In this work we consider matroids to be pairs $M=(E,\mathcal{I})$
	where $E$ is a finite set and $\mathcal{I}$ is the collection of
	independent sets in $M$.  
	We denote by $M/e$ the matroid we get when we contract an
	element $e$ in $M$, and $M \backslash e$ the matroid we get when we
	delete an element $e$ in $M$.
	
	We assume basic familiarity with matroid theory and oriented matroids. The standard references are \cite{oxleybook,bjornerbook}.
	
	\begin{defn}
		Let $M$ be a matroid of rank $r$ and $H$ a hyperplane.
		We say that $H$ is 
		\begin{itemize}
			\item {\em simple}, if there exists $e \in H$ such that $H \backslash \{e\}$ is a flat.     
			\item {\em independent} if it has exactly $r-1$ elements i.e. if it is simple for every element $e \in H$.
			\item {\em multiple} if it is not simple. 
		\end{itemize}
	\end{defn}
	
	After Hansen's counterexample for the natural generalization of Sylvester's theorem to higher dimension, another generalization was proven by him in \cite{hansen65} where he proved that any real point configuration of any dimension has a simple hyperplane but not necessarily an independent one. Later on Murty \cite{mandelthesis} conjectured that Hansen's theorem is also true for oriented matroids i.e. any simple oriented matroid (matroid without loops or parallels) has a simple hyperplane.
	
	As mentioned earlier, the main problem with Hansen's counterexample is that the 6 point configuration does not correspond to a paving matroid since we have three points lying on the same line.
	
	\begin{defn}
		A {\em uniform matroid} with a ground set $E$ of size $n$ and rank $r$ is a matroid where every subset of $E$ of size $r$ is a basis.  A {\em paving matroid} is a matroid in which every circuit has size either $r$ or $r+1$, where $r$ is the rank of the matroid. In a point configuration, this means that no $r-1$ points lie on a same flat of co-dimension 2.
	\end{defn}
	Note that paving matroids are those matroids which are uniform or only slightly deviate from being uniform. Mayhew et al.~\cite{mayhew.et.al} conjecture that almost all matroids are paving. This is mainly why we find paving matroids to be interesting specially since in rank 3, all simple matroids are paving. The following is immediate:
	\begin{prop}
		\label{prop:paving_closed}
		The class of paving matroid is closed under minors.
	\end{prop}
	
	
	Oriented matroids and pseudoline arrangements are strongly connected by the topological representation theorem.
	
	\begin{thm}[The Topological Representation Theorem \cite{folkman1978, winfried02}]
		Any reorientation class of a rank-3 oriented matroid has a representation as a pseudoline arrangement.
		\label{trt}
	\end{thm}
	
	In this representation, the pseudolines in the pseudoline arrangement are the elements of the oriented matroid. Their intersection points are the hyperplanes of the matroids and their simple points of intersection are exactly the independent hyperplanes in the oriented matroid. Therefore, counting simple points in a pseudoline arrangement is actually counting independent hyperplanes in the oriented matroid.
	
	
	
\begin{figure}
		\centering
		\includegraphics[page=4]{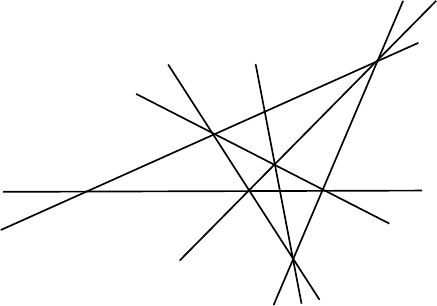}
		\caption{The Kelly-Moser configuration \cite{kellymoser}.}
		\label{kellymoser}
	\end{figure}

	
	Only two main results on the lower bound of simple points in line arrangement were found after Dirac's conjecture in 1951. The first was by Kelly and Moser in 1958 \cite{kellymoser}, where they proved that we have at least 3n/7 simple points in an arrangement of n lines, with equality occurring in the Kelly-Moser configuration shown in Figure \ref{kellymoser}, and the second was by Csima and Sawyer in 1993, \cite{CsimaSawyer1993} where they improved Kelly and Moser's bound to 6n/13 with of course the Kelly-Moser configuration being the only exception. In the same paper, Csima and Sawyer indicated that their proof easily extends to arrangement of pseudolines. This is why we were able to use their result to prove our main result concerning a lower bound for independent hyperplanes in oriented paving matroids.
	
	\begin{thm}\label{thm}
		An oriented paving matroid $M$ of rank $r\geq3$, on n elements, where $n \geq 5+r$,  has at least
		\[ f(n,r)= \dfrac{12}{13(r-1)}\,\,\binom{n}{r-2}\]
		independent hyperplanes.
	\end{thm}
	The proof can be found in Section~\ref{sec:proof}.

	\section{Application of Theorem~\ref{thm}}
	\label{sec:application}
	
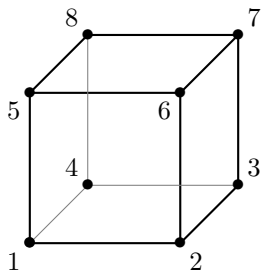
\begin{figure}[h]
			\centering
			\begin{tikzpicture}
				\draw[thick](2,2,0)--(0,2,0)--(0,2,2)--(2,2,2)--(2,2,0)--(2,0,0)--(2,0,2)--(0,0,2)--(0,2,2);
				\draw[thick](2,2,2)--(2,0,2);
				\draw[gray](2,0,0)--(0,0,0)--(0,2,0);
				\draw[gray](0,0,0)--(0,0,2);
				\draw(0,0,0) node[above left] {4};
				\draw(0,0,2) node[below left] {1};
				\draw(2,0,2) node[below right] {2};
				\draw(2,0,0) node[above right] {3};
				\draw(0,2,2) node[below left] {5};
				\draw(2,2,2) node[below left] {6};
				\draw(2,2,0) node[above right] {7};
				\draw(0,2,0) node[above left] {8};
				\foreach \Point in {(0,0,0),(0,0,2),(0,2,0),(0,0,2),(2,0,0),(0,2,2),(2,2,0),(2,0,2),(2,2,2)}{
					\node at \Point {\textbullet};}
			\end{tikzpicture}
			
			\caption{An illustration of the non-orientable simple paving matroid $AG(3,2)'$. (Check \cite{oxleybook} page 508)}
			\label{fig:AG32}
		\end{figure}

		We apply Theorem~\ref{thm} to give a short proof that the Matroid $AG(3,2)'$ is not orientable. The matroid $AG(3,2)'$ shown in Figure~\ref{fig:AG32}, is the unique relaxation of the binary affine cube AG(3,2). It's an 8-element matroid of rank 4 where the 4 point planes are the six faces of the cube, the six diagonals such as \{1,2,7,8\}, and one twisted plane \{1,8,3,6\}.

		$AG(3,2)'$ is a simple paving matroid of rank 4. By our theorem, in order for $AG(3,2)'$ to be orientable, it should have at least
		\[ f(8,4)= \dfrac{12}{13 \cdot 3}\,\,\binom{8}{2}=\dfrac{336}{39} \cong 8.6\]
		independent hyperplanes.
		The independent hyperplanes in this matroid are the sets of 3 elements that do not lie on a common 4-point plane. 
		After looking at the figure we notice that there are only 4 independent hyperplanes, namely: \{2,4,5\}, \{2,4,7\}, \{5,7,2\} and \{5,7,4\}.
		Therefore, by Theorem~\ref{thm}, $AG(3,2)'$ is not orientable. This is already known by da Silva in \cite{dasilva}, where she actually proves a more general and stronger result, but by a more complicated argument.

	\section{Proof of Theorem~\ref{thm}}
	\label{sec:proof}

	We proceed by induction on $r\ge 3$.
	
	By the Topological Representation Theorem~\ref{trt}, a rank-3
	oriented matroid has a representation as a pseudoline arrangement in the projective plane. By Csima and Sawyer \cite{CsimaSawyer1993}, except for the Kelly-Moser configuration, a pseudoline arrangement has a least $6n/13$ simple points. Since $n \geq 5+r$ the Kelly-Moser configuration is excluded. In matroid terms,
	this corresponds to $6n/13$ independent hyperplanes in a rank $3$ oriented matroid founding our induction.
	
	Now assume $r>3$. We fix an element $e \in E$. 
	By Proposition~\ref{prop:paving_closed}, 
	the contraction $M/e$ is paving and has at least 5+(r-1) elements, so by
	induction, $M/e$ has at least 
	\[f(n-1,r-1)=\dfrac{12}{13(r-2)}\,\,\binom{n-1}{r-3}\]
	independent hyperplanes.
	
	Now any independent hyperplane in $M/e$ can be extended by $e$ to an independent hyperplane in $M$. To prove this, we take an independent hyperplane $H$ in $M/e$ and prove that $H\cup e$ is an independent hyperplane in $M$.
	
	We have that $r(H)=r-2$, thus
	$r(H\cup e)=r-1=|H\cup e|$. Therefore $H\cup e$ is independent.
	It is also a flat, because if it wasn't i.e. if there exists an element $x \in \cl(H \cup e) \backslash (H \cup e) $ then $x \in \cl_{M/e}(H)\backslash H$ contradicting,
	$H$ being closed in $M/e$. Therefore $H\cup e$ is an independent flat of rank $r-1$, thus an independent hyperplane in $M$.
	
	Since each of the $f(n-1,r-1)$ independent hyperplanes in $M/e$ can be extended by~$e$ 
	to obtain an independent hyperplane in $M$ and  
	since this is the case for any element $e$ of the matroid, each of these independent hyperplanes will be counted $r-1$ times, and so $M$ has at least
	\[
	\frac{n}{r-1} f(n-1,r-1) = \dfrac{12}{13(r-1)}\,\,\binom{n}{r-2}=f(n,r)
	\] independent hyperplanes.
	This concludes the proof of Theorem~\ref{thm}.

	\section{Discussion} 
	\label{sec:Discussion}
	
	In the rank three case all simple matroids are
	paving, thus being paving is a restriction only in rank at least four.
	While in the pseudoline case the bound is known to be almost
	tight 
	\cite{CsimaSawyer1993}, one might expect, that due to the simplicity of
	the proof our bound should be quite weak for higher dimension. But
	considering $n-1$ points in general position in $\mathbb{R}^{d-1}$ and an
	extra point increasing the dimension, this configuration yields a
	paving matroid of rank $d+1$ with one multiple hyperplane and
	$\binom{n-1}{d-1} = \binom{n-1}{r-2}$ independent ones.
	Note, that
	in this example all but one hyerplane are independent. It still seems
	possible, that always at least a linear fraction of the hyperplanes in
	rank at least $4$ is independent. 
	
	As we have seen, paving matroids are closed under minors. But they are
	not closed under taking duals. A paving matroid where the dual is
	paving as well is called {\em sparse paving}. It is easy to see that a
	paving matroid is sparse paving if for any two circuits $C_1,C_2$ such
	that $|C_1|=|C_2|=r$ we must have $|C_1\cap C_2| \le r-2$. 
	
	One can be interested in the following two problems:
	
	\begin{prob}
		Is there a class of sparse paving matroids of rank $4$ with only a
		subcubic number of independent hyperplanes?
	\end{prob}
	
	Finally, coming back to the simple hyperplanes we would be interested in the following
	
	\begin{prob}
		Does any paving matroid of rank at least $4$ always have $r-2$ points
		which belong to at least as many simple as multiple hyperplanes?
	\end{prob}
	
	Note that in rank 3 the matroid of $K_4$ is a counterexample to
	the above.

\subsection*{Acknowledgment}

The authors would like to thank Manfred Scheucher for his collaboration and efforts which helped to enhance the presentation of the results.
	
	
	
	
	
	\bibliographystyle{plain}
	\bibliography{Independent_hyperplanes_in_oriented_paving_matroids}

\end{document}